\documentclass[11pt,bezier]{article}
\usepackage{amsmath,amssymb,amsfonts,euscript,graphicx}

\newtheorem{The}{Theorem}[section]
\newtheorem{Pro}[The]{Proposition}
\newtheorem{Lem}[The]{Lemma}
\newtheorem{Cor}[The]{Corollary}

\newtheorem{Examp}[The]{Example}

\title{\bf \Large  Modules Satisfying the Prime Radical Condition and a Sheaf Construction
for Modules II\thanks {The research
 of first author   was in part supported by
a grant from IPM (No. 90160034).}
\thanks {2010{ \it Mathematics Subject Classification}: 13C13, 13C99, 13A99, 14A25.
 }}
\author{{\normalsize \bf {M. Aghasi}$^{\rm a}$, { M. Behboodi$^{{\rm a,b}}$\thanks{Corresponding
 author.} and M. Sabzevari$^{{\rm a}}$}}\\
 {\small{ $^{\rm a}$Department of Mathematical Sciences, Isfahan University of Technology}}\vspace{-1mm}\\
  {\small{ P.O.Box: 84156-83111, Isfahan, Iran}}\\
 {\small{ $^{\rm b}$School of Mathematics, Institute for Research in Fundamental Sciences
 (IPM)}}\vspace{-1mm}\\ {\small{ P.O.Box: 19395-5746, Tehran, Iran}}\vspace{-1mm}\\
{\footnotesize{$\mathsf{m.aghasi@cc.iut.ac.ir}$}}\vspace{-1mm}\\
{\footnotesize{$\mathsf{mbehbood@cc.iut.ac.ir}$}}\vspace{-1mm}\\
{\footnotesize{$\mathsf{sabzevari@math.iut.ac.ir}$}}}
\date{}
\begin{document}
\maketitle
\begin{abstract}
{\small \noindent In this paper we continue our study of modules
satisfying the prime radical condition ($\mathbb{P}$-radical
modules), that was introduced in Part I (see \cite{BS}). Let $R$
be a commutative ring with identity. The purpose of this paper is
to   show that the theory of spectrum of $\mathbb{P}$-radical
$R$-modules (with the Zariski topology) resembles to that of
rings. First, we investigate the behavior of $\mathbb{P}$-radical
modules under localization and direct sums. Finally,  we describe
the construction of a structure sheaf on the prime spectrum
Spec$(M)$, which generalizes the classical structure sheaf of the
ring $R$ in Algebraic Geometry to the module $M$. }
\end{abstract}

\section{Introduction}

The present paper is a sequel to \cite{BS} and so the notations
introduced in Introduction of \cite{BS} will remain in force. In
particular, all rings  are associative commutative with identity
$1\neq 0$ and modules are unital. Let $M$ be an $R$-module.  If
$N$ is a submodule (resp. proper submodule) of a module  $M$ we
write $N\leq M$ (resp. $N\lneqq M$). Also, we denote the
annihilator of the factor module $M/N$  by $(N:M)$ and we call $M$
faithful if $(0:M)=0.$

 Recall that the {\it spectrum }
Spec($R$) of a ring $R$ consists of all prime ideals of $R$ and is
non-empty. For each
 ideal $I$ of $R$, we set ${V(I)=\{{\cal{P}}\in {\rm Spec}(R) :
I\subseteq {\cal{P}}\}}$. Then the sets $V(I)$, where $I$ is an
ideal of $R$,  satisfy the axioms for the closed sets of a
topology on Spec($R$), called the {\it Zariski topology} of $R$.
The distinguished open sets of Spec($R$) are the open sets of the
form $D(f)=\{{\cal{P}}\in {\rm Spec}(R) : f\not\in
{\cal{P}}\}={\rm Spec}(R)\setminus V(f)$, where $V(f)=V(Rf)$.
These sets form a basis for the Zariski topology on Spec$(R)$ (see
for examples, Atiyah and Macdonald \cite{Atiyah} and Hartshorne
\cite{Hart}).

A proper submodule $P$ of $M$ is called a {\it prime submodule} if
$am\in P$ for  $a\in R$, and $m\in M$, implies that  $m\in P$ or
$aM\subseteq P$.  The prime spectrum Spec$(M)$ is defined to be
the set of all prime submodules of $M$.

  In the literature, there are many different
generalizations of the Zariski topology of rings to modules via
prime submodules. In  a series of papers (see for example,
\cite{Lu2,Lu3,MMS1,MMS2,MMS3,BH1,BH2,BN}) a group of algebraists
including mainly  Lu,  McCasland,  Moore, Smith and  Behboodi
carried out an intensive and systematic study of the spectrum of
prime submodules. For instance, as [13], for any submodule $N$ of
$M$ we define the {\it variety}
  of $N$, denoted by $V(N)$  to be the set of all prime submodules $P$ of
$M$ such that  $(N:M)\subseteq (P:M)$. Assume that $\zeta(M)$
denotes the collection of all subsets $V(N)$ of Spec$(M)$, then
$\zeta(M)$ contains the empty set and Spec$(M)$, and also,
$\zeta(M)$ is closed under arbitrary intersections and finite
union  and  under finite unions. Thus it is evident that for every
module $M$ there always exists a topology, $\tau$ say, on
Spec$(M)$ having  we conclude that $\zeta(M)$ as the family of all
closed sets. $\tau$ is called the {\it Zariski topology} on
Spec$(M)$. This topology on Spec$(M)$ is studied in
 [13].

It is well-known that for any commutative ring $R$, there is a
sheaf of rings on Spec$(R)$, called the {\it structure sheaf},
denoted by ${\cal{O}}_{{\rm Spec}(R)}$, defined as follows: for
each prime ideal ${\cal{P}}$ of $R$, let $R_{\cal{P}}$ be the
localization of $R$ at ${\cal{P}}$. For an open set $U\subseteq
{\rm Spec}(R)$, we define ${\cal{O}}_{{\rm Spec}(R)}(U)$ to be the
set of functions $s: U\rightarrow \coprod_{{\cal{P}}\in
U}R_{\cal{P}}$, such that $s({\cal{P}})\in R_{\cal{P}}$, for each
${\cal{P}}\in U$, and such that $s$ is locally a quotient of
elements of $R$: to be precise, we require that for each
${\cal{P}}\in U$, there is a neighborhood $V$ of ${\cal{P}}$,
contained in $U$, and elements $a,f\in R$, such that for each
${\cal{Q}}\in  V$, $f\not\in {\cal{Q}}$, and
$s({\cal{Q}})=\frac{a}{f}$ in $R_{\cal{Q}}$ (see for example
Hartshorne \cite{Hart},
 for definition and  basic properties of the sheaf
${\mathcal O}_{{\rm Spec}(R)}$).

For an $R$-module $M$ let $\psi: {\rm Spec}(M)\rightarrow {\rm
Spec}(R/{\rm Ann}(M))$ defined by $\psi(P)=(P:M)/{\rm Ann}(M)$ for
every $P\in {\rm Spec}(M)$. $\psi$ is called the natural map of
${\rm Spec}(M)$. An $R$-module $M$ is called {\it primeful}  if
either $M=(0)$ or $M\neq (0)$ and the natural map $\psi: {\rm
Spec}(M)\rightarrow {\rm Spec}(R/{\rm Ann}(M))$ is surjective.
This notion of primeful module has been extensively studied by Lu
\cite{Lu1}). Let $M$ be a primeful faithful $R$-module.   In
\cite{Tekir}, the author obtained  an $R$-module
$\mathcal{O}_X(U)$ for each open set $U$ in ${\rm Spec}(M)$ such
that $\mathcal{O}_{{\rm Spec}(M)}$ is a sheaf of modules over $X$.
In fact, $\mathcal{O}_{{\rm Spec}(M)}$ is a generalization of the
structure sheaf of rings to primeful faithful modules.

 The purpose of this paper  is to develop the
structure sheaf of  rings to a wider class of modules  called {\it
$\mathbb{P}$-radical modules} that was introduced in Part I (see
\cite{BS}).  Let $N$ be a proper submodule of an $R$-module $M$.
Then the {\it prime radical} $\sqrt[p]{N}$  is the intersection of
all prime submodules of $M$ containing $N$ or, in case there are
no such prime submodules, $\sqrt[p]{N}$  is $M$. Clearly
$V(N)=V(\sqrt[p]{N})$. We note that, for each ideal $I$ of $R$,
$\sqrt[p]{I}=\sqrt{I}$ (the intersection of all prime ideals of
$R$ containing $I$).  Let $M$ be  an $R$-module $M$. Then $M$ is
called  a {\it $\mathbb{P}$-radical module}  if $M$ satisfies
 the prime radical condition
``$(\sqrt[p]{{\cal{P}}M}:M)={\cal{P}}$'' for every prime  ideal
${\cal{P}}\supseteq {\rm Ann}(M)$.   In  \cite[Proposition 2.3 and
Example 2.4]{BS}
 it is shown that the  class of
$\mathbb{P}$-radical  modules contains the class of primeful
$R$-modules properly. Hence, all free modules, all finitely
generated modules and all homogeneous semisimple modules lie in
the class of $\mathbb{P}$-radical modules. Moreover, if $R$ is a
Noetherian ring (or an integral domain), then every projective
$R$-modules is $\mathbb{P}$-radical (see \cite[Theorem 2.5 and
Corollary 2.4(iii)]{BS}). In particular, if $R$ is an Artinian
ring, then
 all $R$-modules
are $\mathbb{P}$-radical (see \cite[Theorem 2.13]{BS}).

In this paper, we continue the study of this construction. In
Section 2, the behavior of $\mathbb{P}$-radical modules under
localization and direct sums  are studied. Section 3 is devoted to
the discussion of the main topic of this paper, the resemblance of
the theory of spectrum of $\mathbb{P}$-radical $R$-modules (with
the Zariski topology) and that of commutative  rings. In fact, for
each $R$-module $M$, we obtain an $R$-module ${\mathcal O}(U)$ for
each open set $U$ in Spec$(M)$. We show that ${\mathcal O}_{{\rm
Spec}(M)}$ is a sheaf of $R$-modules over Spec$(M)$. Moreover, the
following results are two fundamental
theorems in Section 3:\vspace{2mm}\\
\noindent{\bf Theorem 3.1} {\it Let $R$ be a ring and $M$ be an
$R$-module. For any prime submodule $P$ of $M$, the stalk
${{\mathcal O}_{{\rm Spec}(M)}}_P$ is isomorphic to
$M_{\mathcal{P}}$, where
$(P:M)={\mathcal{P}}$.}\vspace{2mm}\\
\noindent{\bf Theorem 3.2} {\it Let $R$ be a ring and  $M$ be a
$\mathbb{P}$-radical   $R$-module. Then for each $f\in R$, the
module ${\mathcal O}_{{\rm Spec}(M)}(D(fM))$ is isomorphic to the
 localized module  $M_f$. In particular, ${\mathcal O}_{{\rm Spec}(M)}({\rm Spec}(M))\cong M$.}

Finally, as an application in Section 4, we prove the following
interesting result:\vspace{2mm}\\
\noindent{\bf Theorem  4.1} {\it Let M be a $\mathbb{P}$-radical
$R$-module and $f,g\in R$. Then, $M_f\cong M_g$ as $R$-modules if
and only if $\sqrt{(fM:M)}=\sqrt{(gM:M)}$.} \

\section{The Behavior of $\mathbb{P}$-Radical Modules Under
Localization and Direct Sums }

 Now we are going to   investigate
  two important properties
of the $\mathbb{P}$-radical   Modules.
 In fact we show the affect of the localization on the $\mathbb{P}$-radical   property of modules.
 Moreover, we investigate the finite direct sum of the $\mathbb{P}$-radical   Modules. First, we  need the following two
lemmas  for our purpose:
\begin{Lem} Let $R$ be a ring and  $I$, $J$ and $K$ be ideals of $R$. Then,
$$\sqrt{(I+J)\cap (I+K)}=\sqrt{I+(J\cap K)}.$$
\end{Lem}
\noindent {Proof.} Let $r\in \sqrt{(I+J)\cap (I+K)}$. Then there
exist two integers $s$ and $t$ such that $r^s=i_1+j$ and
$r^t=i_2+k$ for some $i_1,i_2\in I$, $j\in J$ and $k\in K$. Then
we have $r^{s+t}=i_1i_2+i_1k+i_2j+kj$. So $r^{s+t}\in I+J\cap K$
and consequently $r\in \sqrt{I+(J\cap K)}$ i.e. $\sqrt{(I+J)\cap
(I+K)}\subseteq\sqrt{I+(J\cap K)}$. It is easy to check the
converse of this relation and so the equality satisfied.$~\square$

If $P$ is a prime  submodule of an $R$-module $M$, then
${\cal{P}}= (P:M)$  is a prime ideal and we shall say that $P$  is
a ${\cal{ P}}$-prime submodule.

 \begin{Lem} {\rm \cite [Lemma 4.6]{Lu1}} Let $R$ be a ring, $M_1$ and $M_2$ be $R$-modules  and $\mathcal P\in V(Ann(M))$. If $P$ is a
$\mathcal P$-prime submodule of $M_1$, then $P\bigoplus M_2$
 is a $\mathcal P$-prime submodule of $M$.
\end{Lem}
 \begin{The} Every  finite direct sum of $\mathbb{P}$-radical   modules, is a
$\mathbb{P}$-radical module.
\end{The}
\noindent{\bf Proof.} Let $M_1$ and $M_2$ be two
  $\mathbb{P}$-radical   $R$-modules. It is suffices to show that
  the $R$-module $M=M_1\bigoplus M_2$   is a $\mathbb{P}$-radical   module. Let ${\cal{P}}$
be a prime  ideal of $R$ containing ${\rm Ann}(M)={\rm
Ann}(M_1)\bigcap {\rm Ann}(M_2)$. Clearly
$V({\cal{P}}M_i)=V(({\cal{P}}+{\rm Ann}(M_i))M_i)$
for  $i=1,2$. Now using Lemmas 2.1 and 2.2, we have:\\
 $$(\sqrt[p]{{\cal{P}}M}:M)=\bigcap_{P\in V({\cal{P}}M)}(P:M)\subseteq(\bigcap_{P\in V({\cal{P}}M_1)}(P\oplus M_2:M))\bigcap(\bigcap_{P\in
V({\cal{P}}M_2)}(M_1\oplus P:M))$$
 $$~~~~~~~~~~~~~~~~~~=(\bigcap_{P\in V(({\cal{P}}+{\rm Ann}(M_1))M_1)}(P\oplus M_2:M))\bigcap(\bigcap_{P\in V(({\cal{P}}+{\rm Ann}(M_2))M_2)}(M_1\oplus
P:M))$$
 $$~~~=(\bigcap_{P\in V(({\cal{P}}+{\rm Ann}(M_1))M_1)}(P:M_1))\bigcap(\bigcap_{P\in V((I+{\rm Ann}(M_2))M_2)}(P:M_2))$$
 $$~=(\sqrt[p]{({\cal{P}}+{\rm Ann}(M_1))M_1}:M_1)\bigcap (\sqrt[p]{({\cal{P}}+{\rm Ann}(M_2))M_2}:M_2)$$
 $$=\sqrt{{\cal{P}}+{\rm Ann}(M_1)}\bigcap\sqrt{{\cal{P}}+{\rm Ann}(M_2)}~~~~~~~~~~~~~~~~~~~~~~~~~~~~~~$$
 $$=\sqrt{({\cal{P}}+{\rm Ann}(M_1))\bigcap({\cal{P}}+{\rm Ann}(M_2))}~~~~~~~~~~~~~~~~~~~~~~~~~~~~$$
 $$=\sqrt{{\cal{P}}+({\rm Ann}(M_1\bigcap {\rm Ann}(M_2)))}=\sqrt{{\cal{P}}+{\rm Ann}(M)}={\cal{P}}.~~~~$$
Thus $M$ is a $\mathbb{P}$-radical   module.$~\square$\\

In \cite{Man}, the behavior of prime submodules under
localization, are studied. In this section, $R$ is a commutative
ring and $M$ is an   $R$-module. Suppose $S$ is a multiplicatively
closed subset   of $R$ which contains $1$ but not $0$. Let $M_S$
be the  localization of $M$ at $S$. Let $f: M\rightarrow M_S$ be
the   natural map. For any $R_S$-submodule  $Q$ of $M$, we define
  $Q^c=\{m\in M : f(m)\in Q\}$. The above notations are fixed  for this section.

  \begin{Lem} {\rm \cite[Corollary 2.2]{Man}}  Let $M$ be
an $R$-module and \vspace{2mm}\\
 \indent\indent $A=\big\{P~:~P~is~a~prime~R-submodule~of~M~with~M_S\neq P_S \big\}$\vspace{1mm}\\
\indent \indent $B=\big\{Q~:~Q~is~a~prime~R_S-submodule~of~M_S\big\}.$\vspace{1mm}\\
  Then the map $P\mapsto P_S$ is a bijective   order preserving map
  from $A$ to $B$. Its inverse map is given by  $Q\mapsto Q^c$.
  \end{Lem}

Thus by above lemma we conclude that if  $M$ is  an $R$-module and
$S$ is  a multiplicatively closed subset of $R$, then any prime
submodule of $M_S$ is of the form $P_S$ where $P$ is a prime
submodule of $M$ with $(P:M)\bigcap S=\emptyset$.

   We also need the following two lemmas.

 \begin{Lem} {\rm \cite[Corollary 1]{Lu2}} Let S be a multiplicatively closed subset of a
ring R and M an R-module. If P is a prime submodule of M, then
$(P:M)_S=(P_S:_{R_S}M_S)$
\end{Lem}
 \begin{Lem} Let $M$ be an $R$-module and $S$ be a multiplicatively closed subset
 of $R$. Then, for any ideal $I$ and prime ideal $\mathcal P$ of $R$ we
have,\vspace{2mm}\\
 $~~~~~~~~~~~~~~~~~~~~~~~~~~~~~~~\bigcap_{P\in V(IM)}(P:M)_S=\bigcap_{P_S\in
 V(I_SM_S)}(P_S:_{R_S}M_S)$.
 \end{Lem}
\noindent{\bf Proof.}  Let $P\in V(IM)$. Then, $P$ is a prime
submodule of $M$
 and $(P:M)\supseteq I$. Hence, $P_S$ is either equal to $M_S$ or a prime submodule
 of $M_S$. In the first case, $(P:M)_S=R_S$ and it doesn't affect on the intersections. In the second case,
 $(P:M)_S$ could appear in the right intersection by Lemma 2.5. Conversely, let $P_S\in
 V(I_SM_S)$. By Lemma 2.4, we can assume that
 $P$ is a prime submodule of $M$ with $(P:M)\bigcap S=\emptyset$. Thus
 $((P:M)_S)^c=(P:M)$. On the other hand,
 $((P:M)_S)^c=((P_S:_{R_S}M_S))^c\supseteq (I_S)^c\supseteq I$ and so $(P:M)\supseteq
 I$. i.e $P\in V(IM)$. Hence $(P:M)_S$ will appear
 in the left intersection.$~\square$\\
 \begin{Lem}
Let M be an R-module and S be a multiplicatively closed subset of
R such that for any prime ideal ${\mathcal P}_S\in V(Ann(M_S))$,
$\sqrt[p]{{\mathcal P}_SM_S}\neq M_S$. Then for any    ${\mathcal
P}\in V(Ann(M))$\vspace{2mm}\\
$~~~~~~~~~~~~~~~~~~~~~~~~~~~~~~~~~~~\bigcap_{P\in V({\mathcal
P}M)}(P:M)_S=(\bigcap_{P\in V({\mathcal P}M)}(P:M))_S$.
 \end{Lem}
\noindent{\bf Proof.} First we note  that\vspace{2mm}\\
$~~~~~~~~~~~~~~~~~~\sqrt[p]{{\mathcal P}_SM_S}=\bigcap_{P_S\in
V({\mathcal P_S}M_S)}(P_S:M_S)=\bigcap_{P\in V({\mathcal
P}M)}(P:M)_S\neq
M_S.$\vspace{2mm}\\
 Hence there exists a  prime submodule $P\in V({\mathcal
P}M)$ such that $(P:M)\bigcap S=\emptyset$. Consequently we can
assume that the intersection $\bigcap_{P\in V({\mathcal
P}M)}(P:M)_S$ is constructed from such prime submodules $P\in
V({\mathcal P}M)$ with $(P:M)\bigcap S=\emptyset$. Now let
$\frac{r}{s}\in\bigcap_{P\in V({\mathcal P}M)}(P:M)_S$. Then for
any prime submodule $P$ in the intersection there exists $u\in
(P:M)$ and $t\in S$ such that $\frac{r}{s}=\frac{u}{t}$,    and so
there exists $h\in S$ with $h(rt-su)=0$. Hence since $0,u\in
(P:M)$ then $hrt\in (P:M)$. But $(P:M)\bigcap S=\emptyset$ and so
$r\in (P:M)$. Now since $P$ have been chosen arbitrary from the
intersection, we conclude that $\frac{r}{s}\in (\bigcap_{P\in
V({\mathcal P}M)}(P:M))_S$. Therefore $\bigcap_{P\in V({\mathcal
P}M)}(P:M)_S\subseteq(\bigcap_{P\in V({\mathcal P}M)}(P:M))_S$.
The converse of this    inequality is easy to check and hence the
equality satisfied.$~\square$\\
\begin{Pro}
Let $S$ be a multiplicatively closed subset of $R$ and  $M$ be a
$\mathbb{P}$-radical   $R$-module. If $\sqrt[p]{{\mathcal
P}_SM_S}\neq M_S$ for each  prime ideal ${\mathcal P}_S\in
V(Ann(M_S))$, then $M_S$    is a $\mathbb{P}$-radical
$R_S$-module.
 \end{Pro}
\noindent{\bf Proof.}  Let ${\mathcal P}_S\supseteq {\rm
Ann}(M_S)$ be a prime ideal of $R_S$ where ${\mathcal P}$ is a
prime ideal of $R$. Since ${\mathcal P}_S\supseteq {\rm
Ann}(M_S)\supseteq ({\rm Ann}(M))_S$ and $1\in S$, then ${\mathcal
P}={({\mathcal P}_S)}^c\supseteq (({\rm Ann}(M))_S)^c\supseteq
{\rm Ann}(M)$ where $({{\mathcal P}_S})^c$ is the contraction of
${\mathcal P}_S$. Using Lemmas 2.4 and 2.5, one
 can check that, for any $P\in V({\mathcal P}M)$, if $S\bigcap(P:M)=\emptyset$, then    $P_S\in
V({\mathcal P}_SM_S)$ and if $S\bigcap(P:M)\neq\emptyset$ then
$P_S=M_S$. Hence,\vspace{2mm}\\
 $~~~~~~~~~~~(\sqrt[p]{{\mathcal P}_SM_S)}:_{R_S}M_S)=\bigcap_{P_S\in V({\mathcal P}_SM_S))}(P_S:_{R_S}M_S)$\\
 $~~~~~~~~~~~~~~~~~~~~~~~~~~~~~~~~~~~~~~\subseteq\bigcap_{P\in V({\mathcal P}M)}(P_S:_{R_S}M_S)$\\
 $~~~~~~~~~~~~~~~~~~~~~~~~~~~~~~~~~~~~~~=\bigcap_{P\in  V({\mathcal P}M)}(P:M)_S$\\
 $~~~~~~~~~~~~~~~~~~~~~~~~~~~~~~~~~~~~~~=(\bigcap_{P\in V({\mathcal P}M)}(P:M))_S$\\
 $~~~~~~~~~~~~~~~~~~~~~~~~~~~~~~~~~~~~~~=(\sqrt[p]{{\mathcal P}M}:M)_S.$\vspace{2mm}\\
  Since $M$    is a $\mathbb{P}$-radical module,  $(\sqrt[p]{{\mathcal P}M}:M)=\sqrt{{\mathcal P}}$ and so
  $(\sqrt[p]{{\mathcal P}_SM_S}:_{R_S}M_S)\subseteq
 \sqrt{{\mathcal P}}_S={\mathcal P}_S.$
 This implies that $M_S$    is a $\mathbb{P}$-radical   $R_S$-module.$~\square$\\

 As an important corollary of the above proposition  we have the following result:
 \begin{Cor}
 Let $M$ be a  $\mathbb{P}$-radical $R$-module and $\mathcal P$ be a prime ideal of $R$.
 If  for any prime ideal
 ${\mathcal Q}_{\mathcal P}\in V(Ann(M_{\mathcal P}))$,  $\sqrt[p]{{\mathcal Q}_{\mathcal P}M_{\mathcal P}}\neq M_{\mathcal P}$, then $M_{\mathcal P}$ is a $\mathbb{P}$-radical $R_{\mathcal P}$-module.
  \end{Cor}

 Now we investigate  the converse of the above process. In fact we
 investigate the existence of the    $\mathbb{P}$-radical   property on $M$ with
 the assumption of the existence of this property on its
 localizations.

 \begin{Pro}
Let $M$ be an $R$-module such that for any maximal ideal
${\mathcal M}\in V(Ann(M))$, $M_{\mathcal M}$ is a
$\mathbb{P}$-radical $R_{\cal{M}}$-module and $Ann(M_{\mathcal
M})=(Ann(M))_{\mathcal M}$. Then $M$    is a $\mathbb{P}$-radical
$R$-module.
 \end{Pro}
\noindent{\bf Proof.} Let ${\mathcal P}\supseteq {\rm Ann}(M)$ be
a prime ideal. Then there exists a maximal ideal ${\mathcal M}\in
V({\rm Ann}(M))$ with ${\mathcal P}\subseteq\mathcal M$. Now since
${\rm Ann}(M_{\mathcal M})=({\rm Ann}(M))_{\mathcal M}$ then we
have ${\mathcal P}_{\mathcal M}\supseteq {\rm Ann}(M_{\mathcal
M})$. Hence by the assumption and using Lemma 2.6, the followings
satisfied:\vspace{2mm}\\
 $~~~~~~~~~~~~(\sqrt[p]{{\mathcal P_{\mathcal M}}M_{\mathcal M}}:M_{\mathcal M})={\mathcal P}_{\mathcal M}~~~\Rightarrow~~~\bigcap_{(P_{\mathcal
M}:{M_{\mathcal M})}\supseteq {\mathcal P}_{\mathcal M}}({\mathcal
P}_{\mathcal M}:M_{\mathcal M})={\mathcal P}_{\mathcal M}$\vspace{2mm}\\
 $~~~~~~~~~~~~~~~~~~~~~~~~~~~~~~~~~~~~~~~~~~~~~~~~~\Rightarrow~~~\bigcap_{(P:M)\supseteq\mathcal P}(P:M)_{\mathcal M}={\mathcal P}_{\mathcal M}$\vspace{2mm}\\
 $~~~~~~~~~~~~~~~~~~~~~~~~~~~~~~~~~~~~~~~~~~~~~~~~~\Rightarrow~~~\bigcap_{(P:M)\supseteq\mathcal{P}}(P:M)_{\mathcal M}\subseteq{\mathcal P}_{\mathcal P}.$\vspace{2mm}\\
 It follows that\vspace{2mm}\\
 $~~~~~~~~~~~~~~~~~~~~~~~~~~~~~~~~~~~(\sqrt[p]{{\mathcal P}M}:M)=\bigcap_{(P:M)\supseteq\mathcal
 P}(P:M)$\vspace{1mm}\\
 $~~~~~~~~~~~~~~~~~~~~~~~~~~~~~~~~~~~~~~~~~~~~~~~~~~~~\subseteq ((\bigcap_{(P:M)\supseteq\mathcal P}(P:M))_{\mathcal
 M})^c$\vspace{1mm}\\
 $~~~~~~~~~~~~~~~~~~~~~~~~~~~~~~~~~~~~~~~~~~~~~~~~~~~~\subseteq (\bigcap_{(P:M)\supseteq\mathcal P}(P:M)_{\mathcal
 M})^c$\vspace{1mm}\\
 $~~~~~~~~~~~~~~~~~~~~~~~~~~~~~~~~~~~~~~~~~~~~~~~~~~~~={{\mathcal P}_{\mathcal P}}^c$\vspace{1mm}\\
 $~~~~~~~~~~~~~~~~~~~~~~~~~~~~~~~~~~~~~~~~~~~~~~~~~~~~={\mathcal P}$\vspace{1mm}\\
Therefore, $(\sqrt[p]{{\mathcal P}M}:M)={\mathcal P}$ and so  $M$
is a $\mathbb{P}$-radical $R$-module.$~\square$

 \begin{The} Let $M$ be an $R$-module such that $M_{\cal P}$ is a non-zero $\mathbb{P}$-radical   $R_{\mathcal P}$-module for
 any ${\cal P}\in V(Ann(M))$. Then $M$    is a $\mathbb{P}$-radical $R$-module.
\end{The}
\noindent{\bf Proof.} Let ${\cal{P}}\supseteq {\rm Ann}(M)$ be a
prime deal of $R$. At first we note that  since $M_{\mathcal
P}\neq 0$, ${\rm Ann}(M_{\mathcal P})\subseteq{\mathcal
P}_{\mathcal P}\neq R_{\mathcal P}$.  Thus   by the assumption we
have $(\sqrt[p]{{\cal{P}}_{\mathcal P}M_{\mathcal P}}:M_{\mathcal
P})={\cal{P}}_{\mathcal P}$.  Now by using Lemma 2.6, we
have:\vspace{2mm}\\
$~~~~~~~~~~~~~~~~~~~~~~~~~~~~~~~~~~~(\sqrt[p]{{\cal{P}}M}:M)=\bigcap_{P\in
V({\cal{P}}M)}(P:M)$\vspace{2mm}\\
$~~~~~~~~~~~~~~~~~~~~~~~~~~~~~~~~~~~~~~~~~~~~~~~~~~~~\subseteq((\bigcap_{P\in
V({\cal{P}}M)}(P:M))_{\mathcal P})^c$\vspace{2mm}\\
$~~~~~~~~~~~~~~~~~~~~~~~~~~~~~~~~~~~~~~~~~~~~~~~~~~~~\subseteq(\bigcap_{P\in V({\cal{P}}M)}(P:M)_{\mathcal P})^c$\vspace{2mm}\\
$~~~~~~~~~~~~~~~~~~~~~~~~~~~~~~~~~~~~~~~~~~~~~~~~~~~~=(\bigcap_{P_{\mathcal P}\in V({\cal{P}}_{\mathcal P}M_{\mathcal P})}(P_{\mathcal P}:_{R_{\mathcal P}}M_{\mathcal P}))^c$\vspace{2mm}\\
$~~~~~~~~~~~~~~~~~~~~~~~~~~~~~~~~~~~~~~~~~~~~~~~~~~~~=((\sqrt[p]{{\cal{P}}_{\mathcal
P}M_{\mathcal P}}:M_{\mathcal
P}))^c={{\mathcal P}_{\mathcal P}}^c$\vspace{2mm}\\
$~~~~~~~~~~~~~~~~~~~~~~~~~~~~~~~~~~~~~~~~~~~~~~~~~~~~=\mathcal P.$\vspace{2mm}\\
Therefore, $(\sqrt[p]{{\mathcal P}M}:M)={\mathcal P}$ and so  $M$
is a $\mathbb{P}$-radical $R$-module.$~\square$

\section{ A Structure Sheaf for Modules}

Next we will define an structure sheaf of modules    on the
topological space Spec$(M)$. For each
 prime submodule $P$ of $M$, let ${\mathcal{P}}=(P:M)$
and $M_{\mathcal{P}}$ be the localization of $M$ at
${\mathcal{P}}$. For an open set $U\subseteq {\rm Spec}(M)$, we
define
$${\mathcal O}_{{\rm Spec}(M)}(U)=\{s:U\rightarrow\coprod_{P\in
U}M_{\mathcal{P}}~s.t.~s(P)\in M_{\mathcal{P}}\},$$ where the maps
$s$ are locally constant i.e. for any $P\in U$, there exists an
open neighborhood $V$ of $P$ where $V\subseteq U$ and $s$ is
constant on $V$.  Now it is clear that sums  of such functions
 as well as  the natural functions  $rs$ where  $r\in R$ and
$(rs)(P):=r(s(P))$ are again such. Thus ${\mathcal O}_{{\rm
Spec}(M)}(U)$ is a unitary  $R$-module. If $V\subseteq  U$ are two
open sets, the natural restriction map ${\mathcal O}_{{\rm
Spec}(M)}(U)\longrightarrow {\mathcal O}_{{\rm Spec}(M)}(U)$ is a
homomorphism of $R$-modules. Now, by definition, it is easy to
check that ${\mathcal O}_{{\rm Spec}(M)}$ is a sheaf of
$R$-modules.

Recall that for any sheaf $\mathcal O$ on a topological space $X$
and for any $x\in X$, the stalk of $\mathcal O$ at $x$, denoted by
${\mathcal O}_x$, is represented by all pairs $<U,s>$ where $U$ is
an open subset of $X$ containing $x$ and $s\in {\mathcal O}(U)$.
Two such pairs $<U,s>$ and $<V,t>$ define the same element of
${\mathcal O}_x$ if and only if there is an open neighborhood $W$
of $x$ with $W\subseteq U\cap V$ and $res_{UW}s=res_{VW}t$.
\begin{The} Let $R$ be a ring and $M$ be an
$R$-module. For any prime submodule $P$ of $M$, the stalk
${{\mathcal O}_{Spec(M)}}_P$ is isomorphic to $M_{\mathcal{P}}$,
where $(P:M)={\mathcal{P}}$.
\end{The}
\noindent {Proof.} First we define a homomorphism from ${{\mathcal
O}_{{\rm Spec}(M)}}_P$ to $M_{\mathcal{P}}$, by sending any local
section $s$ in a neighborhood of $P$ to its value $s(P)\in
M_{\mathcal{P}}$.  This gives a well- defined homomorphism
$\varphi$ from ${{\mathcal O}_{{\rm Spec}(M)}}_P$ to
$M_{\mathcal{P}}$. The map $\varphi$ is surjective, because any
element of $M_{\mathcal{P}}$, can be represented as a quotient
$m/f$ with $m\in M$ and $f\in R\setminus {\mathcal{P}}$.  Then
$D(fM)$ will be an open neighborhood of $P$ (since if $fM\subseteq
P$, then $f\in {\mathcal{P}}$, a contradiction),  and $m/f$
defines a section of ${{\mathcal O}_{{\rm Spec}(M)}}$ over $D(fM)$
whose value at $P$ is the given element. To show that $\varphi$ is
injective, let $U$ be a neighborhood of $P$, and let $s,t\in
{\mathcal{O}}(U)$ be
 elements having the same value $s(P)=t(P)$ at $P$. By the
 definition of our structure sheaf and
 shrinking $U$ if necessary, we may assume that $s=m/f$ and $t=n/g$
on $U$, where $m$, $n\in M$, and $f$, $g\in R\setminus
{\mathcal{P}}$. Since  $m/f$ and $n/g$  have the same image in
$M_{\mathcal{P}}$, it follows from the definition of localization
that there is an $h\in R\setminus {\mathcal{P}}$ such that
$h(gm-fn)=0$ in $M$. Therefore $m/f=n/g$ in every local module
$M_{\mathcal{P}'}$ such that $h,~f,~g\in R\setminus
{\mathcal{P}}'$. But the set of such ${\mathcal{P}}'$ is the open
set $D(fM)\cap  D(gM)\cap D(hM)$, which contains $P$. Hence $s=t$
in a whole neighborhood of $P$, so they have the same stalk at
$P$. So  $\varphi$  is an isomorphism, which completes the
proof.\\

Now we are going to prove an important property of the structure
sheaf of    $\mathbb{P}$-radical   modules. This property could
helps us to have a better investigation of the concept of spectrum
of modules in the algebraic geometry point of view. This property
is as follow.
\begin{The}
Let $R$ be a ring and  $M$ be a $\mathbb{P}$-radical   $R$-module.
Then for each $f\in R$, the module ${\mathcal O}_{Spec(M)}(D(fM))$
is isomorphic to the
 localized module  $M_f$. In particular, ${\mathcal O}_{Spec(M)}(Spec(M))\cong M$.
\end{The}
We prove this theorem in several steps. But at first we need the
following lemmas.
\begin{Lem}Suppose that $M$ be a    $\mathbb{P}$-radical   $R$-module. Then for two ideals $I_1$ and $I_2$ of $R$ containing
$Ann(M)$, $V(I_1M)\subseteq V(I_2M)$ implies that
$\sqrt{I_2}\subseteq\sqrt{I_1}$.
\end{Lem}
\noindent {Proof.}  The equality $V(I_1M)\subseteq V(I_2M)$
implies that $\cap_{P\in V(I_1M)}(P:M)\supseteq \cap_{P\in
V(I_2M)}(P:M)$. Thus $Rad_M(I_1)\supseteq Rad_M(I_2)$ and since
$M$    is $\mathbb{P}$-radical   then
$\sqrt{I_1}\supseteq\sqrt{I_2}$.
\begin{Lem}
Let $M$    be a $\mathbb{P}$-radical   $R$-module. Then for two
ideals $I$ and $J$ of $R$, $V(IM)\subseteq V(JM)$ implies that
$\sqrt{(IM:M)}\supseteq\sqrt{(JM:M)}$.
\end{Lem}
\noindent {Proof.} Let $V(IM)\subseteq V(JM)$, for two ideals $I$
and $J$ of $R$. Then $V((IM:M)M)\subseteq V((JM:M)M)$ and so
$\sqrt{(IM:M)}\supseteq\sqrt{(JM:M)}$  by Lemma 3.3.
\begin{Lem}
Let M be a    $\mathbb{P}$-radical   R-module and I be a radical
ideal of R. Then (IM:M)=I if and only if $Ann(M)\subseteq I$.
\end{Lem}
\noindent {Proof.} Let ${\rm Ann}(M)\subseteq I$. Then by
definition of the    $\mathbb{P}$-radical   property we have
$(IM:M)\subseteq(\sqrt[p]{IM}:M)=I$. On the other hand the
converse is obviously satisfied. Thus we have $(IM:M)=I$. The
converse of the statement is natural.
\begin{Lem} Let M be an R-module, $f\in R$ and $s\in{\mathcal O}_{Spec(M)}(D(fM))$. Then there exist some $h_i\in R$
and $m_i\in M$ such that $D(fM)=\bigcup_iD(h_iM)$ and
$s\equiv\frac{m_i}{h_i}$ on each $D(h_iM)$.
\end{Lem}
\noindent {Proof.} Let $s\in {\mathcal O}_{{\rm Spec}(M)}(D(fM))$.
By definition, there exist open subsets $V_i$ of $D(fM)$, on which
$s$ is represented by a quotient $\frac{m_i}{g_i}$, with
$g_i\notin (P:M)$ for all $P\in V_i$. Hence for each $i$,
$V_i\subseteq D(g_iM)$. Now, the open sets of the form $D(hM)$
form a base for the topological space ${\rm Spec}(M)$. Hence we
can assume that $V_i=D(h_iM)$ for some $h_i\in R$. Since
$D(h_iM)\subseteq D(g_iM)$, then we have $D(h_ig_iM)=D(h_iM)\cap
D(g_iM)=D(h_iM)$. On the other hand,
$\frac{m_i}{g_i}=\frac{h_im_i}{h_ig_iM}$. Thus, replacing $h_i$ by
$h_ig_i$ and $m_i$ by $h_im_i$, we can assume that $D(fM)$ is
covered by the open sets $D(h_iM)$ and that $s$ is represented by
$\frac{m_i}{h_i}$ on $D(h_iM)$.
\begin{Lem}Let M be a $\mathbb{P}$-radical   R-module, $f\in R$ and $D(fM)=\bigcup_iD(h_iM)$
for some $h_i\in R$. Then, there exists an integer $n\in\Bbb{N}$
such that $f^n=\sum_{i=1}^kr_ib_i$, for some $r_i\in (h_iM:M)$,
$b_i\in R$ and $k\in\Bbb{N}$.
\end{Lem}
\noindent {Proof.} Let $N:=\sqrt{(\sum(h_iM:M))}M$. Since
$D(fM)\subseteq\cup D(h_iM)$ we have:
\begin{eqnarray*}
&&V((fM:M)M)=V((f)M)\supseteq \cap V((h_i)M)=\\
&&\cap V((h_iM:M)M)=V((\sum(h_iM:M))M)=\\
&&V((\sqrt{\sum(h_iM:M)})M)=V(N)= V((N:M)M).
\end{eqnarray*}
 Hence
$\sqrt{(fM:M)}\subseteq\sqrt{(N:M)}$, by Lemma 3.3. Furthermore we
have $\sqrt{\sum(h_iM:M)}\supseteq$ {\rm Ann}$(M)$ and so
$(N:M)=\sqrt{\sum(h_iM:M)}$, by Lemma 3.5. Hence there exists
$n\in\Bbb{N}$ such that $f^n=b_1r_1+\ldots+b_kr_k$, for some
$b_i\in R$ and $r_i\in (h_iM:M)$.
\begin{Lem}Let M be a $\mathbb{P}$-radical R-module. By the notations of Lemmas
\emph{3.6} and \emph{3.7}, $D(fM)$ can be covered by the open
subsets $D(r_iM)$, $i=1,\ldots,k$. Furthermore, for any
$i=1,\ldots,k$, there exists $n_i\in M$ such that
$s\equiv\frac{n_i}{r_i}$ on each $D(r_iM)$.
\end{Lem}
\noindent {Proof.} Since $r_i\in (h_iM:M)$, for any
$i=1,\ldots,k$, we get $r_iM\subseteq h_iM$. Thus
$V(h_iM)\subseteq V(r_iM)$ and equivalently $D(r_iM)\subseteq
D(h_iM)$. Hence $D(fM)\supseteq \cup_{i=1}^{k} D(r_iM)$. For the
converse, we have to show that $V((f)M)\supseteq \cap_{i=1}^k
V((r_i)M)=V((\sum_{i=1}^k (r_i))M)$. If $P\in V((\sum_{i=1}^k
(r_i))M)$, then $(P:M)\supseteq \sum_{i=1}^k (r_i)$ and in
particular $r_1b_1+\ldots+r_kb_k\in (P:M)$. Hence, by Lemma 3.7,
$f^n\in (P:M)$ and so $P\in V((f^n)M)=V((f)M)$ as wanted.
Furthermore, $r_iM\subseteq h_iM$ implies that there exists
$n_i\in M$ such that $r_im_i=h_in_i$ and so
$\frac{m_i}{h_i}=\frac{n_i}{r_i}$ for $i=1,\ldots,k$. Since
$D(r_iM)\subseteq D(h_iM)$ and $s$ is represented by
$\frac{m_i}{h_i}$ on $D(h_iM)$, we can assume that $s$ is
represented by the quotient $\frac{n_i}{r_i}$ on $D(r_iM)$,
$i=1,\ldots,k$.\\

Now, we are ready to prove the Theorem 3.2. To do this, we define
the canonical map $\psi: M_f\rightarrow {\mathcal O}_{{\rm
Spec}(M)}(D(fM))$ by sending $\frac{m}{f^n}$ to the section $s\in
{\mathcal O}_{{\rm Spec}(M)}(D(fM))$ with $s(P)=\frac{m}{f^n}$,
for any $P\in D(fM)$. Clearly $\psi$ is an $R$-module
homomorphism. By this definition of $\psi$, we complete the proof
of Theorem 3.2 as follows:\vspace{.25cm}\\
\noindent{\bf Proof of Theorem 3.2} It is enough to show that the
homomorphism $\psi$ introduced as above is bijective. At firs we
prove the injectivity of this map. Let
$\psi(\frac{m_1}{f^{n_1}})=\psi(\frac{m_2}{f^{n_2}})$, then for
any $P\in D(fM)$, $\frac{m_1}{f^{n_1}}=\frac{m_2}{f^{n_2}}$ on
$M_{\cal P}$, where ${\mathcal P}=(P:M)$. Thus, for any $P\in
D(fM)$, there exists $h\notin (P:M)$ such that
$h(f^{n_1}m_2-f^{n_2}m_1)=0$. Let $\mathcal{A}$ be the ideal {\rm
Ann}$(f^{n_1}m_2-f^{n_2}m_1)$. Then for any $P\in D(fM)$, there
exists $h\in\mathcal{A}$ such that $h\notin (P:M)$. Thus
$\mathcal{A}\nsubseteq (P:M)$ and so $V((\mathcal{A}M:M)M)\cap
D(fM)=\emptyset$. Furthermore, one can check that
$V((\mathcal{A}M:M)M)=V((\sqrt{\mathcal{A}}M:M)M)$. Consequently,
$V((\sqrt{\mathcal{A}}M:M)M)\subseteq V((f)M)=V((fM:M)M)$. Hence
by Lemma 3.3, $\sqrt{(fM:M)}\subseteq
\sqrt{({\sqrt\mathcal{A}}M:M)}$ and so $f\in
\sqrt{(\sqrt{\mathcal{A}}M:M)}$. On the other hand, by Lemma 3.5,
$(\sqrt{\mathcal{A}}M:M)=\sqrt{\mathcal{A}}$. Thus $f\in
\sqrt{\mathcal{A}}$ and there exists $k\in\Bbb{N}$ with
$f^k(f^{n_1}m_2-f^{n_2}m_1)=0$. Therefore,
$\frac{m_1}{f^{n_1}}=\frac{m_2}{f^{n_2}}$ on $M_f$ and
consequently $\psi$
 is injective. On the other hand, $\psi$ is surjective. Indeed,  According to Proposition 3.7, it is
enough to show that the map $\psi$ considered as above is
surjective. Just as we saw in Lemmas 3.4 to 3.6, there exist
$r_i,b_i\in R$ and $n_i\in M$, $i=1,\ldots,k$, where
$D(fM)=\cup_{i=1}^kD(r_iM)$, $s|_{D(r_iM)}\equiv\frac{n_i}{r_i}$
and $f^n=\sum_{i=1}^kr_ib_i$ for some $n\in\Bbb{N}$. First, we
claim that for any $i,j=1,\ldots,k$, the equality $r_in_j=r_jn_i$
holds. Since $D(r_iM)\cap D(r_jM)=D(r_ir_jM)$, we have two
elements of $M_{r_ir_j}$, namely $\frac{n_i}{r_i}$ and
$\frac{n_j}{r_j}$ both of which represent $s$. Thus, according to
the injectivity of $\psi$, applied to $D(r_ir_jM)$, we have
$\frac{n_i}{r_i}=\frac{n_j}{r_j}$ in $M_{r_ir_j}$. Hence for some
$t\in\Bbb{N}$, the equality $(r_ir_j)^t(r_jn_i-r_in_j)=0$ holds.
Since there are only $k$ indices involved, we may pick $t$ so
large that it works for all $i,j$ at once. Rewrite this equation
as ${r_j}^{t+1}({r_i}^tn_i)-{r_i}^{t+1}({r_j}^tn_j)=0$. Replacing
each $r_i$ by ${r_i}^{t+1}$ and $n_i$ by ${r_i}^tn_i$, we still
have $s$ represented on $D(r_iM)$ by $\frac{n_i}{r_i}$ and
$r_in_j=r_jn_i$ for all $i,j=1,\ldots,k$. Now put
$m=\sum_{i=1}^kb_in_i$, then
$r_jm=\sum_{i=1}^kr_jb_in_i=\sum_{i=1}^kb_ir_in_j=f^nn_j$. This
implies that $\frac{n_j}{r_j}=\frac{m}{f^n}$ on $D(r_jM)$. So
$\psi(\frac{m}{f^n})=s$ everywhere, which shows that $\psi$ is
surjective. $\square$

\section{An application}

 As we saw in $\S1$, the results of this
paper could give the opportunity to consider the concept of
spectrum of modules in algebraic geometry point of view. But some
other interesting and directed results can yield from them. For
instance, the following theorem presents a criterion for the
existence of an isomorphism between two localizations of a
$\mathbb{P}$-radical module:
\begin{The}
Let M be a    $\mathbb{P}$-radical   R-module and $f,g\in R$.
Then, $M_f\cong M_g$ as R-modules if and only if
$\sqrt{(fM:M)}=\sqrt{(gM:M)}$.
\end{The}
\noindent {Proof.} Let $\varphi:M_f\rightarrow M_g$ be an
$R$-module isomorphism and $P\in D(fM)$. We claim that $P$ is an
element of $D(gM)$. Otherwise, if $P\notin D(gM)$ then
$gM_f\subseteq P_f$ and so $\varphi(P_f)\supseteq
g\varphi(M_f)=gM_g=M_g$ which implies that $\varphi(P_f)=M_g$. But
$P\in D(fM)$ and hence $P_f$ is a prime submodule of $M_f$. Thus,
$\varphi(P_f)$ must be a prime submodule of $M_g$ which is a
contradiction. Hence, $P\in D(gM)$ and so $D(fM)\subseteq D(gM)$.
Using the isomorphism $\varphi^{-1}$, we can prove the reverse
inclusion and hence $D(fM)=D(gM)$. Thus,
$\sqrt{(fM:M)}=\sqrt{(gM:M)}$ by
Lemma 3.4.\\
For the converse, let $\sqrt{(fM:M)}=\sqrt{(gM:M)}$. Then
$D(fM)=D(gM)$ and thus $M_f\cong {\mathcal O}_{{\rm
Spec}(M)}(D(fM))={\mathcal O}_{{\rm Spec}(M)}(D(gM))\cong M_g$ by
Theorem 3.2.$~\square$\\

We remark that there exist exemplars of modules such that the
above result is not true in those cases. For example:
\begin{Examp} {\rm
Suppose that $M$ is the Pr\"{u}fer group $\Bbb{Z}(p^{\infty})$ as
a $\Bbb{Z}$-module. $M$ is a torsion divisible module and so it is
primeless by \cite{MMS1}, Lemma 1.3. Hence it    is not a
$\mathbb{P}$-radical   module. Let $q\neq p$ be a prime integer.
Then we have, $\sqrt{(pM:M)}=\sqrt{(qM:M)}=\Bbb{Z}$. But, one can
check that $\frac{1}{q}+\Bbb{Z}$ is an element of rank $q$ in
$M_q$ Although $M_p$ doesn't have any element of this rank. Hence
$M_p\ncong M_q$ as two $\Bbb{Z}$-modules.}
\end{Examp}


\bigskip
\begin{thebibliography}{99}
\bibitem{Atiyah}
\textrm{M. F. Atiyah, I. Macdonald}, \emph{Introduction to
Commutative Algebra}, Addison-Wesley Pub. Co.,
1969.\vspace{-2.5mm}
\bibitem{BH2}
\textrm{M. Behboodi, M. R. Haddadi} \emph{Classical Zariski
topology  of  modules and  spectral spaces I}, Int.  Electron. J.
Algebra {\bf 4} (2008), 104-130.\vspace{-2.5mm}
\bibitem{BH1}
\textrm{M. Behboodi, M. R. Haddadi} \emph{Classical  Zariski
topology of modules and spectral spaces II}, Int.  Electron. J.
Algebra {\bf 4} (2008), 131-1484.\vspace{-2.5mm}
\bibitem{BN}
\textrm{M. Behboodi, M. J. Noori}, \emph{Zariski-like topology on
the classical prime spectrum of a modules}, Bull. Iranian Math.
 Soc. {\bf 35} (1) (2009), 255-271.\vspace{-2.5mm}
\bibitem{BS}
\textrm{M. Behboodi, M. Sabzevari}, \emph{Modules satisfying the
prime radical condition and a sheaf construction for modules I},
Submited.\vspace{-2.5mm}
\bibitem{Cort}
\textrm{ B. Cortzen, L. W. Small}, \emph{Finite extensions of
rings}, Proc. Amer. Math. Soc. {\bf 103} (1988),
1058-1062.\vspace{-2.5mm}
\bibitem{Hart} \textrm{R. Hartshorne}, \emph{Algebraic Geometry},
Springer, New York, 1977.\vspace{-2.5mm}
\bibitem{Lu1}
\textrm{C. P. Lu}, \emph{A module whose prime spectrum has the
suejective natural map}, Houston J. Math, {\bf 33} (1) (2007),
125-143.\vspace{-2.5mm}
\bibitem{Lu2}
\textrm{C. P. Lu}, \emph{Spectra of modules}, Comm. Algebra, {\bf
23} (10) (1995), 3741-3752.\vspace{-2.5mm}
\bibitem{Lu3}
\textrm{C. P. Lu }, {\it The Zariski topology on the prime
spectrum of a module}, Houston J. Math., {\bf 25} (3) (1999),
417-433.\vspace{-2.5mm}
 \bibitem{Man}
 \textrm{S. H. Man}, \emph{One dimensional domains which satisfy the radical formula are Dedekind
domains}, Arch. Math., 66 (1996), 276-279.\vspace{-2.5mm}
\bibitem{MMS1}
\textrm{R. L. McCasland, M. E. Moore, P. F. Smith}, \emph{On the
spectrum of a module over a commutative ring}, Comm. Algebra, {\bf
25} (1) (1997), 79-103.\vspace{-2.5mm}
\bibitem{MMS2}
\textrm{R.L. Mccasland, M.E. Moore, P.F. Smith}, \emph{An
introduction to Zariski spaces over Zariski topologies, Rocky
Mountain J. Math.} {\bf 28} (4)(1998),  1358-1369.\vspace{-2.5mm}
\bibitem{MMS3}
 \textrm{R.L. Mccasland, M.E. Moore,  P.F. Smith},
\emph{Zariski-finite modules}, Rocky Mountain J. Math. {\bf 30}
(2)(2000) 689-701.\vspace{-2.5mm}
 \bibitem{MS}
\textrm{R. L. McCasland, P. F. Smith}, \emph{Prime submodules of
Noetherian modules}, Rocky Mountain J. Math. {\bf 23} (1993),
1041-1062.\vspace{-2.5mm}
\bibitem{Yilmaz}
\textrm{P. F. Smith, D. P. Yilmaz}, \emph{Radicals of submodules
of free modules}, Comm. Algebra, {\bf 27} (5) (1999),
2253-2266.\vspace{-2.5mm}
\bibitem{Tekir}
\textrm{U. Tekir}, \emph{On the sheaf of modules}, Comm. Algebra,
{\bf 33} (2005), 2557-2562.

\end{thebibliography}
\end{document}